\documentclass{amsart}
\renewcommand\newsymbol[5]{%
\DeclareMathSymbol#1{#3}{\ifcase #2\or AMSa\or AMSb\fi}{"#4#5}}
\newcommand\newbbbletter[2]{%
\DeclareMathSymbol#1{0}{AMSb}{`#2}}
\let\emptyset\undefined
\let\le\undefined
\let\ge\undefined
\newsymbol\restriction 1216   
\newsymbol\comb        124F   
\newsymbol\le          1336   
\newsymbol\ge          133E   
\newsymbol\emptyset    203F   
\newbbbletter\reals    R
\newcommand\half{\frac12}
\newcommand\quart{\frac14}
\newcommand\thrquart{\frac34}

\newcommand\preim{^{-1}}
\newcommand\0{\mathbf{0}}                 
\newcommand\1{\mathbf{1}}                 
\newcommand\cl{\operatorname{cl}}         
\newcommand\Int{\operatorname{int}}       
\newcommand\abs[1]{\mathopen|#1\mathclose|}
\theoremstyle{plain}
\newtheorem{thm}{Theorem}[section]
\newtheorem{lem}[thm]{Lemma}
\newtheorem{prop}[thm]{Proposition}
\newcommand\citenames[2]{\textsc{#1} and \textsc{#2}}
\makeatletter
\renewenvironment{cases}%
    {\left\{\,\vcenter\bgroup
              \normalbaselines\openup\jot\m@th
              \ialign\bgroup$##\hfil$&\quad##\hfil\crcr}%
    {\crcr\egroup\egroup\right.}
\makeatother
\begin{document}

\title{Hereditary indecomposability and the Intermediate Value Theorem}

\author{Alan Dow}
\address{Department of Mathematics\\York University\\4700 Keele Street\\
         Toronto, Ontario\\Canada M3J 1P3}
\email{dowa@yorku.ca}

\author[Klaas Pieter Hart]{Klaas Pieter Hart*}
\thanks{*The research of the second author was supported by Netherlands
        Organization for Scientific Research (NWO) --- Grant~R\,61-322.}
\address{Faculty of Information Technology and Systems\\TU Delft\\
         Postbus 5031\\2600~GA~\kern0pt~Delft\\the Netherlands}
\email{k.p.hart@its.tudelft.nl}
\urladdr{http://aw.twi.tudelft.nl/\~{}hart}



\begin{abstract}
We show that hereditarily indecomposable spaces can be characterized
by a special instance of the Intermediate Value Theorem in their ring of
continuous functions.
\end{abstract}

\maketitle

\section*{Introduction}

The classical Intermediate Value Theorem (IVT for short) states that
if $f$~is a continuous function from the interval~$[a,b]$ to~$\reals$ with
$f(a)\cdot f(b)<0$ then there is $c$ in~$(a,b)$ such that~$f(c)=0$.
In \cite{HenriksenLarsonMartinez96}
\citenames{Henriksen, Larson}{Martinez} investigated forms of this theorem
in lattice-ordered rings, where, because of the absence of any natural topology,
they restricted their attention to polynomials.
We mention some of their results for the ring $C^*(X)$ of bounded
real-valued continuous functions on the topological space~$X$;
let us call $X$~an IVT-space if the ring $C^*(X)$ satisfies the Intermediate
Value Theorem (the precise formulation of the IVT in this context follows
below).
The results are:
\begin{enumerate}
\item every IVT-space is an $F$-space;
\item every compact and zero-dimensional $F$-space is an IVT-space;
\item every compact IVT-space is hereditarily indecomposable.
\end{enumerate}
In this note we establish a partial converse to this last result in that
we show that every compact hereditarily indecomposable space satisfies
the IVT for a restricted class of polynomials.

\section{Preliminaries}

We shall only deal with rings of the form $C^*(X)$, so we can, for the time
being restrict our attention to compact Hausdorff spaces.

\subsection{The intermediate value theorem}

In the ring $C^*(X)$ the IVT takes on the following form:
Let $p$~be a polynomial with coefficients in~$C^*(X)$ and let $u$ and $v$
be elements of~$C^*(X)$ such that $p(u)\le\0\le p(v)$, where
$\0$~denotes the zero function.
Then there is $w\in C^*(X)$ such that $u\wedge v\le w\le u\vee v$
and $p(w)=\0$.
The reason for working with $u\wedge v$ and $u\vee v$ is of course that
it is usually not the case that $u(x)\le v(x)$ for all~$x$
(or $v(x)\le u(x)$ for all~$x$).

To get some feeling for what the IVT says in this context let
$p\in C^*(X)[t]$, so $p(t)=\sum_{i=0}^nf_i t^i$ for some elements~$f_0$,
\dots,~$f_n$ of~$C^*(X)$, and let $u,v\in C^*(X)$ be such that
$p(u)\le \0\le p(v)$.
For every separate~$x\in X$ we get an ordinary polynomial
$p_x(t)=\sum_{i=0}^nf_i(x) t^i$; and the assumptions on $u$ and~$v$
imply that $p\bigl(u(x)\bigr)\le 0\le p\bigl(v(x)\bigr)$.
The classical IVT therefore guarantees that there is a function $w:X\to\reals$
such that $u\wedge v\le w\le u\vee v$ and $p(w)=\0$;
the IVT for $C^*(X)$ demands that this $w$ be continuous.

That this puts severe restrictions on the space~$X$ may be seen as follows:
let $f\in C^*(X)$ and consider the polynomial $p(t)=\abs{f}t-f$.
Now $p(\1)=\abs{f}-f\ge0$ and $p(-\1)=-\abs{f}-f\le\0$, so if $X$ is
an IVT-space there must be a continuous function~$w$ such that
$-\1\le w\le\1$ and $f=w\abs{f}$.
This however is one of the characterizations of $F$-spaces ---
see \citenames{Gillman}{Jerison} \cite{GillmanJerison60}.

\subsection{Hereditarily indecomposable spaces}

Much of what follows is taken from \citenames{Oversteegen}{Tymchatyn}
\cite{OversteegenTymchatyn86}, which is a convenient survey on
hereditarily indecomposable spaces.

To begin we recall that a continuum is said to be \emph{indecomposable} if
it cannot be written as the union of two proper subcontinua; it is
\emph{hereditarily indecomposable} if every subcontinuum is indecomposable.

We use the following characterization of hereditarily indecomposable
continua.

\begin{thm}
A continuum $X$ is hereditarily indecomposable if and only if whenever two
disjoint closed sets $A$ and $B$ and open neighbourhoods $U$ and~$V$
respectively are given we can write $X$ as the union of three closed sets
$X_0$, $X_1$ and~$X_2$ such that
$A\subseteq X_0$,
\ $B\subseteq X_2$,
\ $X_0\cap X_1\subseteq V$,
\ $X_0\cap X_2=\emptyset$, and
  $X_1\cap X_2\subseteq U$.
\end{thm}

The property in this theorem can also be used to characterize those
compact spaces (connected or not) for which every closed connected
subspace is indecomposable; we shall call these compact spaces
hereditarily indecomposable as well.
Observe that with this definition compact zero-dimensional spaces are
hereditarily indecomposable as well.

\section{The IVT implies hereditary indecomposability}
\label{sec.IVT.implies.her-ind}

In this section we reprove Theorem~3.2 from
\citenames{Henriksen, Larson}{Martinez} \cite{HenriksenLarsonMartinez96},
which states that compact IVT-spaces are hereditarily indecomposable.
In their proof these authors used a polynomial of degree~$7$ with
two potentially irreducible quadratic factors.
We use a completely factored polynomial of degree~$3$.

\begin{thm}
Compact IVT-spaces are hereditarily indecomposable.
\end{thm}

\begin{proof}
Let $X$ be a compact IVT-space.
To show that $X$ is hereditarily indecomposable we take disjoint closed
sets $A$~and~$B$ and open sets $U$ and~$V$ such that $A\subseteq U$ and
$B\subseteq V$.
We must exhibit three closed sets $X_0$, $X_1$ and~$X_2$ such that
$A\subseteq X_0$, \ $B\subseteq X_2$, \ $X_0\cap X_1\subseteq V$,
\ $X_0\cap X_2=\emptyset$, \ $X_1\cap X_2\subseteq U$
and $X_0\cup X_1\cup X_2=X$.

Choose a continuous function $f:X\to[0,1]$ such that
$f\restriction A\equiv0$,
\ $f\restriction B\equiv1$,
$f\preim\bigl[[0,\half)\bigr]\subseteq U$
and $f\preim\bigl[(\half,1]\bigr]\subseteq V$.
Using~$f$ we define three continuous functions, $f_1$, $f_2$ and~$f_3$, as
follows:
first
$$
f_1(x)=\begin{cases}
       f(x)-\quart   & if $f(x)\le\quart$\cr
       0             & if $\quart\le f(x)\le\thrquart$\cr
       f(x)-\thrquart& if $\thrquart\le f(x)$;\cr
       \end{cases}
$$
second
$$
f_2(x)=\begin{cases}
       \half\bigl(f(x)-\quart\bigr)   & if $f(x)\le\quart$\cr
       2\bigl(f(x)-\quart\bigr)       & if $\quart\le f(x)\le\thrquart$\cr
       \half\bigl(f(x)+\frac54\bigr)& if $\thrquart\le f(x)$;\cr
       \end{cases}
$$
and third
$$
f_3(x)=\begin{cases}
       f(x)+\thrquart& if $f(x)\le\quart$\cr
       1             & if $\quart\le f(x)\le\thrquart$\cr
       f(x)+\quart   & if $\thrquart\le f(x)$;\cr
       \end{cases}
$$
(At this point the reader may find it instructive to draw the graphs of
 $f_1$, $f_2$ and~$f_3$ in case $X=[0,1]$ and $f(x)=x$.
 The zig-zag that appears when one follows the graph of~$f_3$ left-to-right
 until it meets the graph of~$f_2$ then follows the graph of~$f_2$
 right-to-left until it meets the graph of~$f_1$ and finally the graph of~$f_1$
 left-to-right until the end is characteristic of hereditarily indecomposable
 spaces.)
Note that $f_1\le f_2\le f_3$.

Consider the polynomial $p$ defined by $p(t)=(t-f_1)(t-f_2)(t-f_3)$.
Then $p(\0)\le\0\le p(\1)$, for one readily checks that
\begin{itemize}
\item $f_2(x)<0<f_3(x)<1$ if $f(x)<\quart$;
\item $f_1(x)=0$ and $f_3(x)=1$ if $\quart\le f(x)\le\thrquart$ and
\item $0<f_1(x)<1<f_2(x)$ if $f(x)>\thrquart$.
\end{itemize}
An application of the Intermediate Value Theorem gives us a continuous
function $w:X\to[0,1]$ such that $p(w)=\0$.

Let $X_0=\bigl\{x:w(x)=f_3(x)\bigr\}$,
\ $X_1=\{x:w(x)=f_2(x)\bigr\}$ and $X_2=\bigl\{x:w(x)=f_1(x)\bigr\}$.
We check that these sets have all the required properties.
\begin{itemize}
\item $X_0\cup X_1\cup X_2=X$ because $p(w)=\0$;
\item $A\subseteq X_0$ because if $x\in A$ then $f(x)=0$, hence $w(x)=f_3(x)$;
\item $B\subseteq X_2$ because if $x\in B$ then $f(x)=1$, hence $w(x)=f_1(x)$;
\item $X_0\cap X_1\subseteq V$ because if $x\in X_0\cap X_1$ then
      $f_3(x)=w(x)=f_2(x)$ hence $f(x)=\thrquart$ and $x\in V$;
\item $X_1\cap X_2\subseteq U$ because if $x\in X_1\cap X_2$ then
      $f_1(x)=w(x)=f_2(x)$ hence $f(x)=\quart$ and $x\in U$ and
\item $X_0\cap X_2=\emptyset$ because $f_3-f_1=\1$.
\end{itemize}
We conclude that $X$ is indeed hereditarily indecomposable.
\end{proof}

As announced before, in the next section we shall see that hereditary
indecomposability is in fact characterized by the particular instance
of the Intermediate Value Theorem that was actually employed.

\section{Hereditary indecomposability implies part of the IVT}
\label{sec.her-ind.implies.IVT}

In this section we show that every compact hereditarily indecomposable
$F$-space~$X$ satisfies the Intermediate Value Theorem for
\emph{completely factored polynomials},
that is, polynomials that can be written as
$\prod_{i=1}^n(t-f_i)$, where the $f_i$ are elements of~$C(X)$.

This is a rather limited class of polynomials of course but, as we saw
in Section~\ref{sec.IVT.implies.her-ind}, the case~$n=3$ is already strong
enough to imply hereditary indecomposability.
The Intermediate Value Theorem for this class of polynomials therefore
characterizes hereditary indecomposability for $F$-spaces.

So let $X$ be a hereditarily indecomposable $F$-space and
let $p$, defined by~$p(t)=\prod_{i=1}^n(t-f_i)$, be a completely
factored polynomial in~$C(X)$.
Assume furthermore that $u,v\in C(X)$ are such that $p(u)\le\0\le p(v)$.
Through a series of reductions we show that there is $w\in C(X)$
such that $p(w)=\0$ and $u\wedge v\le w\le u\vee v$.

\begin{lem}
We may assume that $f_1\le f_2\le\cdots\le f_n$.
\end{lem}

\begin{proof}
For each $i\le n$ define $g_i$ by
$$
g_i=\bigwedge_{|F|=i}\bigvee_{j\in F}f_j.
$$
Observe that $g_1\le g_2\le\cdots\le g_n$ and that, for each individual~$x$,
the sets of values $\bigl\{g_1(x),g_2(x),\ldots,g_n(x)\bigr\}$ and
$\bigl\{f_1(x),f_2(x),\ldots,f_n(x)\bigr\}$ are equal.
It follows from this that the coefficients of~$t^0$, $t^1$, \dots,~$t^{n-1}$
in $\prod_{i=1}^n(t-f_i)$ and $\prod_{i=1}^n(t-g_i)$ are the same
and hence that the polynomials are the same.
\end{proof}

The case $n=1$ should offer no problems and the case $n=2$ is dealt with
in the following proposition, which is a special case of
Theorem~2.3\thinspace(b) of
\citenames{Henriksen, Larson}{Martinez} \cite{HenriksenLarsonMartinez96}.
In fact, the polynomial~$p$ need not even be factored; it can always
be factored by completing the square.

\begin{prop}\label{prop.quadratic}
Every space satisfies the Intermediate Value Theorem for monic quadratic
polynomials.
\end{prop}

\begin{proof}
Let $p(t)=t^2+2ft+g$ be such a polynomial and assume that there are~$u$ and~$v$
such that $p(u)\le\0\le p(v)$.
Completing the square gives us $q(t)=(t+f)^2+g-f^2$.
Now because $p(u)\le\0\le p(v)$ we know that $f^2-g\ge\0$ so that we can write
$f^2-g=h^2$ for some nonnegative $h\in C(X)$.
We find that $p(t)=(t+f-h)(t+f+h)$; write $-f-h=f_1$ and $-f+h=f_2$.

Observe that for each~$x$ either $v(x)\le f_1(x)$ or $v(x)\ge f_2(x)$
and that $f_1(x)\le u(x)\le f_2(x)$.
We cover our space by three closed sets:
$P=\cl\bigl\{x:u(x)<v(x)\bigr\}$,
\ $Q=\bigl\{x:u(x)=v(x)\bigr\}$ and
$R=\cl\bigl\{x:u(x)>v(x)\bigr\}$.
We now note that $u\le f_2\le v$ on~$P$
(because $u(x)\le f_2(x)\le v(x)$ whenever $u(x)<v(x)$)
and that $v\le f_1\le u$ on~$R$.
We define $w$ as the combination
$$
(f_2\restriction P)\comb (u\restriction Q)\comb(f_1\restriction R).
$$
Note that $w$ is well-defined because, by continuity, $u\equiv v\equiv f_2$
on~$P\cap Q$ and $u\equiv v\equiv f_1$ on~$Q\cap R$.
Also $p(w)(x)=0$ for all~$x$; this is clear on~$P\cup R$ and on~$Q$ it holds
because $p(u)(x)\le0\le p(v)(x)$ and $p(u)(x)=p(v)(x)$.
Finally, $w$~is continuous because it is the combination of continuous
functions defined on closed subsets.
\end{proof}

We have given such an extensive proof of Proposition~\ref{prop.quadratic}
because it contains elements that we will use quite often in what follows,
to wit breaking the space into closed pieces according to the position
of the~$f_i(x)$ with respect to $u(x)$ and~$v(x)$, and defining $w$ by cases.
From now on we assume that $n\ge3$.

To begin, for every~$x$ we have $f_n(x)\ge u(x)\ge f_{n-1}(x)$
or $f_{n-2}(x)\ge u(x)\ge f_{n-3}(x)$ etc., \emph{because} $p(u)(x)\le0$;
if, for example, $f_{n-1}(x)>u(x)>f_{n-2}(x)$ then clearly $p(u)(x)>0$.
This sequence ends with $f_2(x)\ge u(x)\ge f_1(x)$ if $n$~is even
and with $f_1(x)\ge u(x)$ if $n$~is odd.

Likewise, for all $x$ we have $v(x)\ge f_n(x)$ or
$f_{n-1}(x)\ge v(x)\ge f_{n-2}(x)$ or \dots~or
$f_1(x)\ge v(x)$ if $n$~is even and $f_2(x)\ge v(x)\ge f_1(x)$ if $n$~is odd.

We shall also employ the cover of $X$ by the sets
$P=\cl\bigl\{x:u(x)<v(x)\bigr\}$,
\ $Q=\bigl\{x:u(x)=v(x)\bigr\}$ and
$R=\cl\bigl\{x:u(x)>v(x)\bigr\}$.
On~$Q$ there is no choice: the only admissible solution is
$w_Q=u\restriction Q=v\restriction Q$.
However, once we have found solutions~$w_P$ on~$P$ and $w_R$ on~$R$
then $w=w_P\comb w_Q\comb w_R$ is the desired solution.
On $P$ we have $u\le w_P\le v$ so by continuity we know that
$u(x)=w_P(x)=v(x)$ for all $x\in P\cap Q$.
Likewise $u(x)=w_R(x)=v(x)$ for all $x\in Q\cap R$.
Thus, $w$~is well-defined and as a combination of continuous functions defined
on closed subsets it is continuous.

Because hereditary indecomposability is a closed hereditary property
we can work inside $P$ and $R$ respectively without worrying about the rest
of~$X$.

\subsection{Reduction to odd~$n$}

Assume $n$~is even and recall that in this case $f_1\le u\le f_n$.

We show that on~$P$ we have $q(u)\le0\le q(v)$, where
$q(t)=\prod_{i=2}^n(t-f_i)$.
Indeed the possible positions of~$u(x)$ ensure that $q(u)(x)\le0$ for
all~$x$.
Also, for all~$x$ with $u(x)<v(x)$ we have $v(x)\ge f_2(x)$ because
$f_1(x)\le u(x)<v(x)<f_2(x)$ would imply $p(v)(x)<0$.
Hence, by continuity, $v\ge f_2$ on~$P$, so that $q(v)\ge\0$ on~$P$.

On the set $R$ we can show in a similar fashion that $v\le f_{n-1}$
and hence that $r(u)\ge\0\ge r(v)$, where $r(t)=\prod_{i=1}^{n-1}(t-f_i)$.

Both $q$ and $r$ are of degree~$n-1$.

\medskip
From now on we assume $n\ge 3$ and $n$~odd.

\subsection{Reduction to $u\le v$}

If $u(x)>v(x)$ then, because $f_n\ge u$, we must have $v(x)\le f_{n-1}(x)$
and because $v\ge f_1$ we must have~$u(x)\ge f_2(x)$.
So on~$R$ we get, by continuity, $v\le f_{n-1}$ and $u\ge f_2$.
Consider now the polynomial $q(t)=\prod_{i=2}^{n-1}(t-f_i)$.
Because of the possible positions for $u(x)$ and $v(x)$ listed above
we conclude that $q(u)\ge\0\ge q(v)$.

\subsection{The final case}

We now show how to produce $w$, given that
\ 1)~$X$ is the closure of~$\bigl\{x:u(x)<v(x)\bigr\}$,
\ 2)~$p(u)\le p(v)$ and
\ 3)~$n$~is odd.

Let $k$ be such that $n=2k+1$.
For each $i\le k$ consider the closed sets
$A_i=\cl\bigl\{x:v(x)<f_{2i+1}(x)\bigr\}$
and $B_i=\cl\bigl\{x:u(x)>f_{2i-1}(x)\bigr\}$.

Note that, because of the positioning of the values $u(x)$ and $v(x)$
we have $A_i\subseteq C_i=\bigl\{x:v(x)\le f_{2i}(x)\bigr\}$
and $B_i\subseteq D_i=\cl\bigl\{x:u(x)\ge f_{2i}(x)\bigr\}$.
Now note that
$C_i\cap D_i\subseteq \bigl\{x:u(x)=f_{2i}(x)=v(x)\bigr\}$;
as the set on the right-hand side is nowhere dense
it follows that $\Int C_i$ and $\Int D_i$ are disjoint.

Also, because $X$~is an $F$-space, we know that
$A_i\subseteq\Int C_i$ and $B_i\subseteq \Int D_i$.

Now apply hereditary indecomposability to find three closed sets $X_i$,
$Y_i$ and $Z_i$ that cover~$X$ and with the following properties:
$A_i\subseteq X_i$,
\ $B_i\subseteq Z_i$,
\ $X_i\cap Y_i\subseteq \Int D_i$,
\ $Y_i\cap Z_i\subseteq \Int C_i$
and $X_i\cap Z_i=\emptyset$.
We note the following facts:
\begin{enumerate}
\item $u\le f_{2i-1}$ on $X_i\cup Y_i$ because this set is disjoint
      from~$B_i$;\label{i}
\item $v\ge f_{2i+1}$ on $Y_i\cup Z_i$ because this set is disjoint
      from~$A_i$;\label{ii}
\item $u=f_{2i-1}=f_{2i}$ on $X_i\cap Y_i$ because this set is contained
      in~$D_i$ and because of~(\ref{i});\label{iii}
\item $v=f_{2i+1}=f_{2i}$ on $Y_i\cap Z_i$ because the set is contained
      in~$C_i$ and because of~(\ref{ii}); and\label{iv}
\item $u\le f_{2i-1}\le f_{2i}\le f_{2i+1}\le v$ on~$Y_i$ because
       of (\ref{i})~and~(\ref{ii}).\label{v}
\end{enumerate}
Now we are ready to define~$w$.
We start by letting $w=f_1$ on~$X_1$ and $w=f_2$ on~$Y_1$.
We continue by letting, for $i>1$,
\ $w=f_{2i-1}$ on~$X_i\cap\bigcap_{j<i}Z_j$
and $w=f_{2i}$ on~$Y_i\cap\bigcap_{j<i}Z_j$.
Finally, on $\bigcap_{i\le k}Z_i$ we let~$w=f_n$.

\subsubsection*{We check that $w$ is well-defined.}
By~(\ref{iii}) we have $f_{2i-1}=f_{2i}$ on~$X_i\cap Y_i$ for every~$i$.
If $j<i$ then $X_j\cap X_i\cap\bigcap_{l<i}Z_l=\emptyset$;
on $Y_j\cap X_i\cap\bigcap_{l<i}Z_l\subseteq Y_j\cap Z_j\cap Z_{i-1}$ we have
$v=f_{2j+1}=f_{2j}$ and $v\ge f_{2i-1}$ and so $f_{2i-1}=f_{2j}$,
and on $Y_j\cap Y_i\cap\bigcap_{l<i}Z_l\subseteq Y_j\cap Z_j\cap Z_{i-1}$
we have $v=f_{2j=1}=f_{2j}$ and $v\ge f_{2i+1}\ge f_{2i}$ and
so $f_{2j}=f_{2i}$.
Finally, on $Y_j\cap\bigcap_{i\le k}Z_i\subseteq Y_j\cap Z_j$
we have $v\ge f_n$ and $v=f_{2j}$ so $f_{2j}=f_n$.

\subsubsection*{We check that $u\le w\le v$.}
On $X_1$ we surely have $u\le f_1\le v$ and if $i>1$ then on
$X_i\cap\bigcap_{j<i}Z_j$ we have $u\le f_{2i-1}$ because of~(\ref{i})
and $f_{2i-1}\le v$ because of~(\ref{i}) for~$i-1$.
On each~$Y_i$ we have $u\le f_{2i}\le v$ by~(\ref{v}).
Finally, on $\bigcap_{i\le k}Z_i$ we have $v\ge f_n\ge u$ by~(\ref{ii}).

We see that $w$ is a well-defined continuous function on~$X$ such that
$u\le w\le v$ and, because for all~$x$ there is an~$i$ with $w(x)=f_i(x)$,
such that $p(w)=\0$.

\section{Questions and Conjectures}
\label{sec.q.and.c}

The basic question as to what actually characterizes IVT-spaces remains.
On the basis of the evidence from
Section~\ref{sec.her-ind.implies.IVT}
we conjecture that hereditarily indecomposable spaces also satisfy the
IVT for monic polynomials.

The general case seems more complicated in that the leading coefficient
(and others) may vanish at certain places.
It may very well be that that the full IVT characterizes zero-dimensionality.


\providecommand{\bysame}{\leavevmode\hbox to3em{\hrulefill}\thinspace}

\end{document}